\documentclass[a4paper,12pt]{article}
\usepackage[utf8]{inputenc}
\usepackage{amsfonts, amsopn, amsthm,amssymb,amsmath}
\usepackage{graphicx}
\usepackage{xspace}
\usepackage{hyperref}

\usepackage{tikz,tkz-graph}
\usepackage{subcaption}
\tikzset{vertex/.style={inner sep=1.8pt, outer sep=0pt, circle, fill=black},
emvertex/.style={inner sep=1.8pt, outer sep=0pt, circle, draw, fill=white}}

\newcommand{\G}{\mathcal G}
\newcommand{\Csixv}{C_{6,v}}
\newcommand{\Csixc}{C_{6,c}}
\newcommand{\conf}[1]{Conf(#1)}
\newtheorem{theorem}{Theorem}
\newtheorem{lemma}[theorem]{Lemma}
\newtheorem{corollary}[theorem]{Corollary}
\newtheorem{definition}[theorem]{Definition}

\title{Triangle-free planar graphs with the smallest independence number\thanks{Work on this paper was supported by project 14-19503S (Graph coloring and structure) of Czech Science Foundation.}}
\author{Zdeněk Dvořák\thanks{Computer Science Institute (CSI), Charles University in Prague.
E-mail: \protect\href{mailto:rakdver@iuuk.mff.cuni.cz}{\protect\nolinkurl{rakdver@iuuk.mff.cuni.cz}}.}\and
Tomáš Masařík\thanks{Department of Applied Mathematics, Charles University in Prague.
E-mail: \protect\href{mailto:tarken@kam.mff.cuni.cz}{\protect\nolinkurl{tarken@kam.mff.cuni.cz}}.}\and
Jan Musílek\thanks{Department of Applied Mathematics, Charles University in Prague.
E-mail: \protect\href{mailto:stinovlas@kam.mff.cuni.cz}{\protect\nolinkurl{stinovlas@kam.mff.cuni.cz}}.}\and
Ondřej Pangrác\thanks{Computer Science Institute (CSI), Charles University in Prague.
E-mail: \protect\href{mailto:pangrac@iuuk.mff.cuni.cz}{\protect\nolinkurl{pangrac@iuuk.mff.cuni.cz}}.}}

\begin{document}
\maketitle

\begin{abstract}
Steinberg and Tovey~\cite{SteinbergTovey1993} proved that every $n$-vertex planar triangle-free graph has an independent
set of size at least $(n+1)/3$, and described an infinite class of tight examples.  We show that all $n$-vertex planar triangle-free graphs
except for this one infinite class have independent sets of size at least $(n+2)/3$.
\end{abstract}

By a well-known theorem of Gr\"otzsch~\cite{grotzsch1959}, every planar triangle-free graph is $3$-colorable.
This clearly implies that such a graph $G$ with $n$ vertices has an independent set of size at least $n/3$,
i.e., $\alpha(G)\ge n/3$ in the usual notation.
This can be slightly improved---using a strengthening of Gr\"otzsch's theorem, Steinberg and Tovey~\cite{SteinbergTovey1993}
proved that the equality is never achieved in this bound.
\begin{theorem}[Steinberg and Tovey~\cite{SteinbergTovey1993}]\label{thm-sttov}
If $G$ is an $n$-vertex planar triangle-free graph, then $\alpha(G)\ge (n+1)/3$.
\end{theorem}
They also described an infinite
class $\G$ of planar triangle-free graphs (see Definition~\ref{def:constructionG} below) such that $\alpha(G)=(|V(G)|+1)/3$ for all $G\in\G$.
In this paper, we give a new proof of their result, which also implies that $\G$ contains all the graphs for that the bound is tight (throughout the paper,
we only consider simple graphs without loops or parallel edges).
\begin{theorem}\label{thm-main}
If $G$ is a planar triangle-free graph with $n$ vertices and $G\not\in \G$, then $\alpha(G)\ge (n+2)/3$.
\end{theorem}

Let us mention several related results.  A better known (and much harder) relative of our problem concerns independent
sets in unconstrained planar graphs.  By Four Color Theorem~\cite{AppHak1,AppHakKoc}, each $n$-vertex planar graph has
an independent set of size at least $n/4$.  This bound is tight, and unlike our case, the (infinitely many) known examples
do not seem to exhibit an easily discernible structure.  Indeed, even the algorithmic problem of testing whether
an $n$-vertex planar graph has an independent set greater than $n/4$ has no known polynomial-time solution~\cite{Niedermeier2006,FellowsEtAl2012}.

The \emph{fractional chromatic number} $\chi_f$ of a graph $G$ is the minimum value of $a/b$ over
all positive integers $a\ge b$ for which there exists a coloring that assigns each vertex of $G$ a subset of $\{1,\ldots, a\}$ of size $b$
such that the sets assigned to adjacent vertices are disjoint.  It is easy to see that $\chi_f(G)\le \chi(G)$ and $\alpha(G)\ge |V(G)|/\chi_f(G)$.
Hence, the results above indicate that the fractional chromatic number of $n$-vertex planar triangle-free graphs might be bounded by $3-3/(n+1)$.
As Dvořák et al.~\cite{frpltr} proved, this is the case for planar triangle-free graphs of maximum degree at most $4$; in general, they were
only able to obtain a weaker upper bound $3-3/(3n+1)$.

It is natural to ask whether the bound from Theorem~\ref{thm-main} can be improved, at the expense of having further
families of exceptional graphs.  Algorithmically, this question was answered by Dvořák and Mnich~\cite{dmnich,dmnich-full},
who proved that if an $n$-vertex planar triangle-free graph does not have an independent set larger than $(n+k)/3$, then
its tree-width is $O(\sqrt{k})$.  Using their techique, a more detailed answer can be given,
showing that all such graphs are created from graphs of bounded size by a construction similar to the one used to
define the class $\G$ below; we will give details in a followup paper.  For small values of $k$, an exact description
of exceptional graphs can be obtained using the argument of the current paper (we decided not to present them here, since
the number of exceptional classes grows quickly and dealing with them would obscure the idea).

\section{The extremal graphs}

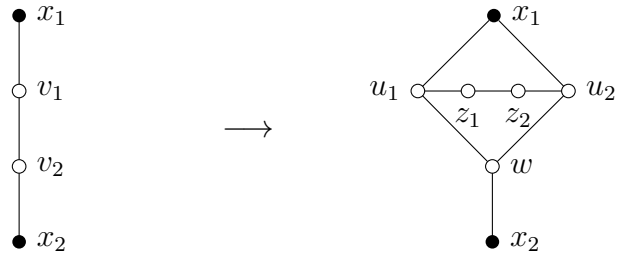
\begin{figure}
\centering
\begin{subfigure}{0.2\textwidth}
\centering
\begin{tikzpicture}[scale=1]
	\draw (0,0) node[vertex,label=right:$x_1$] (x1) {}
	-- ++(-90:1)    node[emvertex,label=right:$v_1$] (v1) {}
	-- ++(-90:1)    node[emvertex,label=right:$v_2$] (v2) {}
	-- ++(-90:1)    node[vertex,label=right:$x_2$] (x2) {};
\end{tikzpicture}
\end{subfigure}
\qquad$\longrightarrow$\qquad
\begin{subfigure}{0.3\textwidth}
\centering
\begin{tikzpicture}[scale=1]
	\draw (0,0) node[vertex,label=right:$x_1$] (x1) {}
	++(-90:1)++(180:1) node[emvertex,label=left:$u_1$] (u1) {}
	++(0:0.66) node[emvertex,label=below:$z_1$] (z1) {}
	++(0:0.66) node[emvertex,label=below:$z_2$] (z2) {}
	++(0:0.66) node[emvertex,label=right:$u_2$] (u2) {}
	++(180:1)++(-90:1)    node[emvertex,label=right:$w$] (w) {}
	-- ++(-90:1)    node[vertex,label=right:$x_2$] (x2) {};
	\draw (u1) -- (z1) -- (z2) -- (u2) -- (x1) -- (u1) -- (w) -- (u2);
\end{tikzpicture}
\end{subfigure}
\caption{Path--diamond replacement.\label{fig-repl}}
\end{figure}

The class $\G$ is defined via the following construction, see Figure~\ref{fig-repl}.  
A $5$-cycle $C=u_1z_1z_2u_2w$ in a graph $G$, where $u_1$ and $u_2$ have a common neighbor $x_1\not\in V(C)$,
$w$ is adjacent to another vertex $x_2\not\in V(C)$, and $\deg_G(u_1)=\deg_G(u_2)=\deg_G(w)=3$ and $\deg_G(z_1)=\deg_G(z_2)=2$, is called a \emph{diamond}.
Let $G_1$ be a graph and let $x_1v_1v_2x_2$ be a path in $G_1$ with $\deg_{G_1}(v_1)=\deg_{G_1}(v_2)=2$.  Let $G_2$ be the graph obtained from the disjoint union of $G_1-\{v_1,v_2\}$ and a $5$-cycle
$u_1z_1z_2u_2w$ by adding the edges $x_1u_1$, $x_1u_2$ and $x_2w$.  We say that $G_2$ is obtained from $G_1$ by a \emph{path--diamond} replacement.  
Conversely, let $G'_2$ be a graph containing a diamond $C=u_1z_1z_2u_2w$, and let $G'_1$ be the graph obtained from $G'_2-V(C)$ by
adding a path $x_1v_1v_2x_2$ with new vertices $v_1$ and $v_2$.  We say that $G'_1$ is obtained from $G'_2$ by
\emph{replacing a diamond by the path $x_1v_1v_2x_2$}.  Note that both of these operations preserve planarity and do not create triangles.

\begin{definition}
\label{def:constructionG}
The class $\G$ consists of the path $P_2$ on two vertices, the $5$-cycle, and all
graphs obtained from the $5$-cycle by a repeated application of the path--diamond replacement.
\end{definition}

Let $C_5^\dagger$ denote the graph obtained from the $5$-cycle by the path--diamond replacement (see Figure~\ref{fig-cdag}),
and let $C_5^\ddag$ denote the graph obtained from $C_5^\dagger$ by the path--diamond replacement (see Figure~\ref{fig-cddag}).
Note that these graphs and their plane drawings are unique up to isomorphism.

\section{Reducible configurations}

An $n$-vertex graph $G$ is \emph{tight} if $G$ is planar, triangle-free, and $\alpha(G)\le(n+1)/3$.

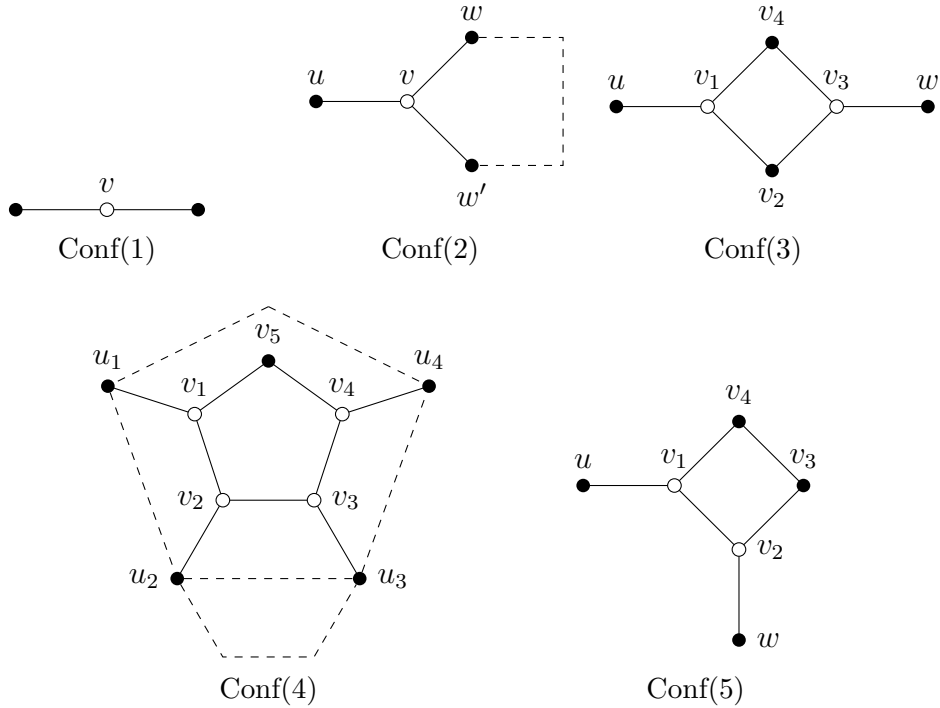
\begin{figure}
\begin{center}
\subcaptionbox{\conf{1}}[0.3\textwidth]{%
\begin{tikzpicture}[scale=1.2]
	\draw (0,0) node[vertex] (u) {}
	-- ++(0:1)    node[emvertex,label=above:$v$] (v) {}
	-- ++(0:1)    node[vertex] (w) {};
\end{tikzpicture}}
\subcaptionbox{\conf{2}}[0.3\textwidth]{%
\begin{tikzpicture}[scale=1.2]
	\draw (0,0) node[vertex,label=above:$u$] (u) {}
	-- ++(0:1) node[emvertex,label=above:$v$] (v) {}
	-- ++(45:1) node[vertex,label=above:$w$] (w) {};
	\draw (v) -- ++(-45:1) node[vertex,label=below:$w'$] (w') {};
	\draw [dashed] (w) -- ++(0:1) -- ++(-90:1.414) -- (w');
\end{tikzpicture}}
\subcaptionbox{\conf{3}}[0.3\textwidth]{%
\begin{tikzpicture}[scale=1.2]
	\draw (0,0) node[vertex,label=above:$u$] (u) {}
	-- ++(0:1) node[emvertex,label=above:$v_1$] (v1) {}
	-- ++(-45:1) node[vertex,label=below:$v_2$] (v2) {}
	-- ++(45:1) node[emvertex,label=above:$v_3$] (v3) {}
	-- ++(135:1) node[vertex,label=above:$v_4$] (v4) {}
	-- (v1);
	\draw (v3) -- ++(0:1) node[vertex,label=above:$w$] (w) {};
\end{tikzpicture}}
\end{center}

\begin{center}
\subcaptionbox{\conf{4}}[0.5\textwidth]{%
\begin{tikzpicture}[scale=1.2]
	\draw (0,0) node[emvertex,label=right:$v_3$] (v3) {}
	-- ++(72:1) node[emvertex,label=above:$v_4$] (v4) {}
	-- ++(144:1) node[vertex,label=above:$v_5$] (v5) {}
	-- ++(216:1) node[emvertex,label=above:$v_1$] (v1) {}
	-- ++(288:1) node[emvertex,label=left:$v_2$] (v2) {}
	-- (v3);

	\draw (v3) -- ++(300:1) node[vertex,label=right:$u_3$] (u3) {};
	\draw (v4) -- ++(18:1) node[vertex,label=above:$u_4$] (u4) {};
	\draw (v1) -- ++(162:1) node[vertex,label=above:$u_1$] (u1) {};
	\draw (v2) -- ++(240:1) node[vertex,label=left:$u_2$] (u2) {};

	\draw [dashed] (v5)++(90:0.6) -- (u4);
	\draw [dashed] (v5)++(90:0.6) -- (u1);
	\draw [dashed] (u2) -- (u3);
	\draw [dashed] (u1) -- (u2);
	\draw [dashed] (u3) -- (u4);

	\draw [dashed] (u2) -- ++(300:1) -- ++(0:1) -- (u3);
\end{tikzpicture}}
\subcaptionbox{\conf{5}}[0.3\textwidth]{%
\begin{tikzpicture}[scale=1.2]
	\draw (0,0) node[vertex,label=above:$u$] (u) {}
	-- ++(0:1) node[emvertex,label=above:$v_1$] (v1) {}
	-- ++(-45:1) node[emvertex,label=right:$v_2$] (v2) {}
	-- ++(45:1) node[vertex,label=above:$v_3$] (v3) {}
	-- ++(135:1) node[vertex,label=above:$v_4$] (v4) {}
	-- (v1);
	\draw (v2) -- ++(270:1) node[vertex,label=right:$w$] (w) {};
\end{tikzpicture}}
\end{center}

\caption{Reducible configurations.\label{fig-redu}}
\end{figure}

We now describe several \emph{reducible configurations} (see Figure~\ref{fig-redu}), which allow local transformations in plane triangle-free graphs that preserve tightness;
this will lead to a natural inductive proof of Theorem~\ref{thm-main}.  For each configuration (with the exception of \conf{5},
which is handled separately in Lemma~\ref{lemma-red5}) we specify the corresponding local transformation, resulting in a \emph{reduced graph}.  We also introduce the notion of
\emph{interference with the outer face}, which is needed later in the proof that one of these configurations appears
in each plane triangle-free graph.

Let $G$ be a plane triangle-free graph with the outer face bounded by a cycle $K$.
\begin{itemize}
\item Configuration \conf{1} consists of a vertex $v\in V(G)$ of degree at most $2$.  The reduced graph is
	obtained by deleting $v$ and all its neighbors.
        The configuration interferes with the outer face if $v\in V(K)$.

\item Configuration \conf{2} consists of a vertex $v\in V(G)$ of degree $3$ with neighbors $u$, $w$, and $w'$,
        such that $G$ contains no path of length $3$ between $w$ and $w'$.  The reduced graph is
	obtained by deleting $u$ and $v$, and by identifying $w$ and $w'$ to a single vertex and suppressing the parallel edges.
	The configuration interferes with the outer face if $\{v,w,w'\}\cap V(K)\neq\emptyset$.

\item Configuration \conf{3} consists of a $4$-face $C=v_1v_2v_3v_4$ in $G$ such that $\deg(v_1)=\deg(v_3)=3$.
	The reduced graph is obtained by deleting $V(C)$ and the neighbors of $v_1$ and $v_3$.
	The configuration interferes with the outer face if $\{v_1,v_3\}\cap V(K)\neq\emptyset$.

\item Configuration \conf{4} consists of a $5$-face $C=v_1v_2v_3v_4v_5$ in $G$ with $\deg(v_1)=\ldots=\deg(v_4)=3$,
        such that, denoting for $i=1,\ldots, 4$ the neighbor of $v_i$ outside of $C$ by $u_i$, $G-V(C)$ contains
	no path of length at most $2$ between $u_1$ and $u_4$, and no path of length $1$ or $3$ between $u_2$ and $u_3$,
	and $u_1u_2,u_3u_4\not\in E(G)$.
	The reduced graph is obtained by deleting $V(C)$, adding the edge $u_1u_4$, and identifying $u_2$ with $u_3$
	to a single vertex and suppressing the parallel edges.
	The configuration interferes with the outer face if $\{v_1,\ldots,v_4,u_1,\ldots,u_4\}\cap V(K)\neq\emptyset$.

\item Configuration \conf{5} consists of a $4$-face $v_1v_2v_3v_4$ in $G$ such that $\deg(v_1)=\deg(v_2)=3$.
	The configuration interferes with the outer face if $\{v_1,v_2\}\cap V(K)\neq\emptyset$.
\end{itemize}

Configuration \conf{5} is dealt with using the following observation.

\begin{lemma}\label{lemma-red5}
Let $G$ be a plane triangle-free graph.  If $G$ contains the configuration \conf{5},
then it also contains the configuration \conf{2}.
\end{lemma}
\begin{proof}
Let $v_1v_2v_3v_4$ be a $4$-face in $G$ with $\deg(v_1)=\deg(v_2)=3$.
If $G$ contains no path of length three between $v_1$ and $v_3$,
or no path of length three between $v_2$ and $v_4$, then \conf{2}
appears in $G$.  However, both such paths cannot be present, since $G$ is plane
and triangle-free.
\end{proof}

Let us now argue that the described reductions preserve tightness.

\begin{lemma}\label{lm:reduction_preserve}
Let $G$ be a plane triangle-free graph containing one of the reducible configurations
\conf{1}, \ldots, \conf{4}, and let $G'$ be the corresponding reduced graph.
Then $G'$ is planar and triangle-free.  Moreover,
there exists a positive integer $k$ such that $|V(G')| \ge |V(G)| - 3k$ and $\alpha(G)\ge\alpha(G') + k$. 
\end{lemma}

\begin{proof}
Let us consider each of the configurations separately; we use the same labels for the vertices
of the configurations as in their definition.  Let $S$ denote the largest independent set in $G'$.

\begin{enumerate}
\item[\conf{1}] We delete $v$ and its (at most two) neighbors, and thus $|V(G')| \ge |V(G)| - 3$.
        Furthermore, $S\cup\{v\}$ is an independent set in $G$, and thus $\alpha(G)\ge\alpha(G') + 1$.
\item[\conf{2}] The identification of $w$ with $w'$ does not create any triangles, since $G$ contains
        no path of length $3$ between these two vertices.  Note that $|V(G')|=|V(G)| - 3$.
	Let $z$ denote the vertex created by the identification of $w$ and $w'$.  If $z\in S$,
	then $(S\setminus\{z\})\cup \{w,w'\}$ is an independent set in $G$; otherwise, $S\cup \{v\}$
	is an independent set in $G$.  Consequently, $\alpha(G)\ge\alpha(G') + 1$.
\item[\conf{3}] Note that $|V(G')|\ge |V(G)| - 6$, and $S\cup \{v_1,v_3\}$ is an independent set
        in $G$, implying $\alpha(G)\ge\alpha(G') + 2$.
\item[\conf{4}] Suppose $G'$ contains a triangle.  Since the distance in $G-V(C)$ between $u_1$ and $u_4$
       is greater than $2$ and $G-V(C)$ contains no path of length $3$ between $u_2$ and $u_3$, we conclude that
       the triangle contains both the edge $u_1u_4$ and the vertex $z$ created by the identification of $u_2$ and $u_3$.
       By planarity, it follows that $u_1u_2,u_3u_4\in E(G)$.  However, this is forbidden by the assumptions
       of the configuration.
       
       Note that $|V(G')|\ge |V(G)| - 6$.  Since $u_1u_4\in E(G')$, by symmetry we can assume
       that $u_1\not\in S$.  If $z\in S$, then $(S\setminus\{z\})\cup\{v_1,u_3,u_4\}$ is an independent set in $G$;
       otherwise, $S\cup \{v_1,v_3\}$ is an independent set in $G$.
       Hence, $\alpha(G)\ge\alpha(G') + 2$.
\end{enumerate}
\end{proof}

\begin{corollary}\label{cor:tight_reduction}
If $G$ is a tight graph containing one of the reducible configurations
\conf{1}, \ldots, \conf{4}, then the corresponding reduced graph $G'$ is also tight.
\end{corollary}
\begin{proof}
By Lemma~\ref{lm:reduction_preserve}, $G'$ is planar and triangle-free.  Furthermore,
there exists $k>0$ such that $|V(G')| \ge |V(G)| - 3k$ and $\alpha(G)\ge\alpha(G') + k$.
Since $G$ is tight, we have $\alpha(G)\le (|V(G)|+1)/3$.  It follows that
$$\alpha(G')\le\alpha(G) - k \le (|V(G)|-3k+1)/3\le (|V(G')|+1)/3.$$
Therefore, $G'$ is also tight.
\end{proof}

\section{Excluding the configurations}

In this section, we argue that tight graphs cannot contain the reducible configurations.
Let us start with some observations on diamonds.

\begin{lemma}\label{lm:reduced_iset}
Let $G$ be a graph containing a diamond $C=u_1z_1z_2u_2w$, and let $G'$ be obtained from $G$ by replacing
the diamond $C$ by the path $x_1v_1v_2x_2$.  Then $|V(G)|=|V(G')|+3$ and $\alpha(G)=\alpha(G')+1$.  Moreover,
for every independent set $S'$ of $G'$, there exists an independent set $S$ of $G$ such that $|S|=|S'|+1$
and $S\setminus V(C)=S'\setminus\{v_1,v_2\}$.
\end{lemma}
\begin{proof}
Consider any independent set $S'$ of $G'$; the independent set $S$ in $G$ of size $|S'|+1$ can be obtained from $S'\cup \{z_2\}$ by
replacing $v_1$ by $u_1$ and replacing $v_2$ by $w$.  Hence, $\alpha(G)\ge\alpha(G')+1$.

Conversely, consider any maximal independent set $S$ of $G$;
note that if $u_1,u_2\in S$, then $S\setminus\{u_2\}\cup \{z_2\}$ is an independent set of the same size, and
if $|\{u_1,u_2\}\cap S|\le 1$, then either $z_1$ or $z_2$ belongs to $S$ by the maximality of $S$.  Hence, by
symmetry we can assume that $z_2\in S$, and an independent set in $G'$ of size $|S|-1$ can be obtained from $S\setminus\{z_2\}$
by replacing $u_1$ by $v_1$ and replacing $w$ by $v_2$.  This implies that $\alpha(G')\ge \alpha(G)-1$.

Combining the inequalities, we conclude that $\alpha(G)=\alpha(G')+1$.
\end{proof}

Let us remark that Lemma~\ref{lm:reduced_iset} implies that an $n$-vertex graph $G\in\G$ satisfies $\alpha(G)=(n+1)/3$.
We say that $G$ is a \emph{minimum counterexample} (to Theorem~\ref{thm-main}) if $G$ is a tight graph not belonging to $\G$ with the smallest
number of vertices (our aim is to prove that no such counterexample exists).

\begin{corollary}\label{cor:missing_conf}
Minimum counterexamples do not contain diamonds.
\end{corollary}
\begin{proof}
Suppose that a minimum counterexample $G$ contains a diamond.  Let $G'$ be the graph obtained from $G$ by replacing the
diamond by a path.  Since $G$ is tight, Lemma~\ref{lm:reduced_iset} implies that $G'$ is tight, and by the minimality of $G$,
we conclude that $G'\in \G$.  However, $G$ is obtained from $G'$ by a path--diamond replacement, and thus $G\in\G$, which is a contradiction.
\end{proof}

Next, we show a useful fact about maximum independent sets in graphs from the class $\G$.

\begin{lemma}\label{lm:avoiding_set}
Consider any graph $G \in \G$, and let $f$ be a face of a plane drawing of $G$ such that $f$ is not incident
with any vertex of degree at most two.  Then there exists an independent set $S\subseteq V(G)$
such that $|S| = (|V(G)|+1)/3$ and $S\cap V(f) = \emptyset$.
\end{lemma}
\begin{proof}
We proceed by the induction on the number of vertices of $G$; hence, assume that the claim holds for
all graphs with less than $|V(G)|$ vertices.
Since $f$ is not incident with any vertices of degree at most two,
$G$ is not $P_2$, the $5$-cycle, or the graph $C_5^\dagger$.

\begin{figure}
\centering
\begin{tikzpicture}[scale=1]
	\draw (0,0) node[emvertex,label=above:$v_1$] (x1) {}
	++(-90:1)++(180:1) node[emvertex,label=left:$v_2$] (u1) {}
	++(0:0.66) node[vertex] (z1) {}
	++(0:0.66) node[emvertex] (z2) {}
	++(0:0.66) node[vertex] (u2) {}
	++(180:1)++(-90:1)    node[emvertex,label=above:$v_3$] (w) {}
	-- ++(-90:1)    node[emvertex,label=below:$v_4$] (x2) {};
	\draw (u2)++(0:1) node[vertex] (a1) {};
	\draw (u1)++(180:1) node[emvertex,label=left:$v_5$] (a2) {}
	++(280:2.6) node[vertex] (a3) {};
	\draw (a1)++(260:2.6) node[emvertex] (a4) {};
	\draw (u1) -- (z1) -- (z2) -- (u2) -- (x1) -- (u1) -- (w) -- (u2);
	\draw (a2) -- (x1) -- (a1) -- (x2) -- (a2) -- (a3) -- (a4) -- (a1);
\end{tikzpicture}
\caption{A maximum independent set in $C_5^\ddag$.}\label{fig-cddag}
\end{figure}
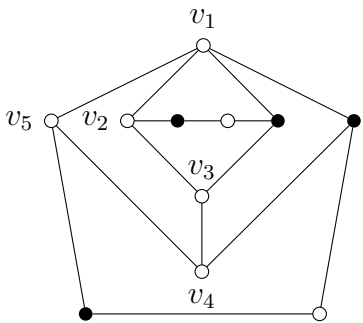

Suppose that $G$ is the graph $C_5^\ddag$.
This graph has a unique plane drawing, with two faces not incident with degree $2$ vertices.
An independent set of size $4$ disjoint from one such face $v_1\ldots v_5$ is depicted in
Figure~\ref{fig-cddag}; the case of the other face is symmetric.

Finally, suppose that $G$ is any other graph in $\G$.  Observe that $G$ contains a diamond $C$ such that
$f$ is not incident with any of the edges of the diamond.  Let $G'$ be the graph obtained from $G$ by
replacing the diamond $C$ by a path, with the natural drawing in the plane preserving the
face $f$.  By the induction hypothesis, $G'$ contains an independent set $S'$ of size $(|V(G')|+1)/3=(|V(G)|+1)/3-1$
disjoint from $V(f)$.  Lemma~\ref{lm:reduced_iset} implies that $G$ contains an independent set of size
$|S'|+1=(|V(G)|+1)/3$ disjoint from $V(f)$.
\end{proof}

We are now ready to show that minimum counterexamples cannot contain \conf{1}.

\begin{lemma}\label{lm:degree2}
A minimum counterexample has minimum degree at least three.
\end{lemma}
\begin{proof}
Suppose that $G$ is a minimum counterexample containing a vertex $v$ of degree $d\le 2$, i.e.,
the configuration \conf{1}.  Let $G'$ be the corresponding reduced graph (obtained from $G$ by removing $v$ and its neighbors), and note that
$|V(G')|=|V(G)|-d-1$.  By Lemma~\ref{lm:reduction_preserve}, we have $\alpha(G)\ge \alpha(G')+1$,
which by Theorem~\ref{thm-sttov} implies $\alpha(G)\ge (|V(G')|+4)/3=(|V(G)|+3-d)/3$.
Since $G$ is tight, we conclude that $d=2$ and $G'$ is tight.  By the minimality of $G$, it follows that $G'\in\G$.

Let $z_1$ and $z_2$ be the neighbors of $v$ in $G$.  Note that there exists a face $f$ of $G'$ such that
the path $z_1vz_2$ of $G$ is drawn within $f$.  Let $N$ denote the set of vertices in $V(f)$ that are adjacent
in $G$ with $z_1$ or $z_2$.  Observe that every maximum independent set of $G'$ intersects $N$, as otherwise
this independent set together with $\{z_1,z_2\}$ would give an independent set in $G$ of size greater than $(|V(G)|+1)/3$.

If $G'$ is a path on two vertices, it follows that $N=V(G')$, and since $G$ is triangle-free, we conclude that $G$ is a $5$-cycle and
$G\in \G$.  If $G'$ is a $5$-cycle $x_1\ldots x_5$, then by symmetry we can assume that $\{x_1,x_2,x_3\}\subseteq N$, and since $G$
is triangle-free, it follows that say $z_1$ is adjacent to $x_1$ and $x_3$, and $z_2$ is adjacent to $x_2$; consequently, $G$ is isomorphic
to $C_5^\dagger$ and $G\in \G$.  In both cases, we obtain a contradiction.

\begin{figure}
\centering
\begin{tikzpicture}[scale=1]
	\draw (0,0) node[vertex,label=above:$v_1$] (x1) {}
	++(-90:1)++(180:1) node[vertex,label=above:$v_2$] (u1) {}
	++(0:0.66) node[vertex,label=above:$v_3$] (z1) {}
	++(0:0.66) node[vertex,label=above:$v_4$] (z2) {}
	++(0:0.66) node[vertex,label=above:$v_5$] (u2) {}
	++(180:1)++(-90:1)    node[vertex,label=below:$u_3$] (w) {}
	-- ++(170:1.5)    node[vertex,label=below:$u_2$] (x2) {};
	\draw (x1)++(190:1.5) node[vertex,label=above:$u_1$] (a2) {};
	\draw (u1) -- (z1) -- (z2) -- (u2) -- (x1) -- (u1) -- (w) -- (u2);
	\draw (x1) -- (a2) -- (x2);
\end{tikzpicture}
\caption{The graph ${C_5^\dagger}$.}\label{fig-cdag}
\end{figure}
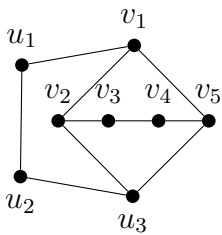

Next, consider the case $G'$ is isomorphic to $C_5^\dagger$; we label its vertices as in Figure~\ref{fig-cdag}.
By symmetry, we can assume that $f=v_1v_2v_3v_4v_5$.  Since $N$ intersects the maximum independent sets $\{u_1, u_3, v_i\}$ for $i\in \{3,4\}$,
and $\{u_1,v_2,v_5\}$, we can also assume that $\{v_2,v_3,v_4\}\subseteq N$.  Since $G$ is triangle-free, it follows that say $z_1$ is
adjacent to $v_2$ and $v_4$, and $z_2$ is adjacent to $v_3$; consequently, $G$ is isomorphic to $C_5^\ddag$ and $G\in \G$, which is again a contradiction.

Finally, consider the case that $G'$ is any other graph in $\G$.  Since $N$ intersects all maximum independent sets of $G'$,
Lemma~\ref{lm:avoiding_set} implies that $f$ is incident with a vertex of degree two of $G'$.  However, then $G'$ contains
a diamond $C$ such that the faces incident with $C$ are distinct from $f$.  It follows that $C$ is also a diamond in $G$,
which contradicts Corollary~\ref{cor:missing_conf}.
\end{proof}

Finally, let us exclude all other configurations.

\begin{lemma}\label{lemma-noredu}
A minimum counterexample does not contain any of the reducible configurations \conf{1}, \ldots, \conf{5}.
\end{lemma}
\begin{proof}
By Lemmas~\ref{lemma-red5} and \ref{lm:degree2}, it suffices to show that a minimum counterexample $G$ does not contain any of
the configurations \conf{2}, \ldots, \conf{4}.  Suppose for a contradiction that $G$ contains one of these configurations,
and let $G'$ be the corresponding reduced graph.  By Corollary~\ref{cor:tight_reduction}, $G'$ is tight, and by the minimality
of $G$, we conclude that $G'\in\G$.  Observe that since $G$ has minimum degree at least three, each reduction results in a graph
with at least two non-adjacent vertices; consequently, $G'\neq P_2$.

Let $D$ denote the set of vertices of $G'$ of degree at most $2$.
Since $G$ has minimum degree at least three, all vertices of degree at most two in $G'$ arise in the reduction.
Observe that either all vertices of $D$ are incident with one face of $G'$ (when \conf{3} is being reduced),
or there exists a vertex $x\in V(G')$ such that all vertices of $D$ that are neither equal nor adjacent to $x$
are incident with one face of $G'$ and form an independent set (when \conf{2} or \conf{4} is being reduced; $x$ is the vertex
created by the identification of two vertices of $G$).  The latter condition is false for all graphs in $\G\setminus\{P_2\}$.
The only graph in $\G\setminus\{P_2\}$ satisfying the former condition is the $5$-cycle.

However, a straightforward case analysis shows that no triangle-free graph of minimum degree at least three containing \conf{3}
reduces to a $5$-cycle.  This is a contradiction.
\end{proof}

\section{Unavoidability}

We finish the proof by showing that every plane triangle-free graph contains one of the configurations \conf{1}, \ldots, \conf{5}.
Our proof is motivated by a similar argument of~\cite{DvoKawTho}.  To deal with short separating cycles, we need to prove a stronger claim.

We say that a vertex in a plane graph $G$ is \emph{internal} if it is not incident with the outer face of $G$.
For a cycle $C$ in a plane graph $G$, let $G_C$ denote the subgraph of $G$ drawn in the closed disk bounded by $C$.
Let $\Csixc$ denote the plane graph consisting of a $6$-cycle that forms its outer face and a chord separating its interior to two $4$-faces, and
let $\Csixv$ denote the plane graph consisting of a $6$-cycle $C$ that forms its outer face and a vertex $v$ adjacent to every other vertex of $C$; see Figure~\ref{fig-C6spec}.
A cycle $C$ in a plane graph $G$ is \emph{dangerous} if its length is at most $6$, $C$ does not bound the outer face of $G$, and $G_C$ is distinct from $C$ itself, $\Csixc$ and $\Csixv$.

\begin{lemma}\label{lemma-unav}
	Let $G$ be a plane triangle-free graph with the outer face bounded by a $(\le\!6)$-cycle $K$, such that $G$ is distinct from $K$ itself, $\Csixc$ and $\Csixv$.
	If $G$ does not contain any dangerous cycle, then it contains one of the configurations \conf{1}, \ldots, \conf{5} that does not interfere with the outer face.
\end{lemma}
\begin{proof}
Suppose for a contradiction that every configuration \conf{1}, \ldots, \conf{5} in $G$ interferes with its outer face.
In particular, since $G$ does not contain \conf{1} not interfering with the outer face, all internal vertices of $G$ have degree at least three (and since $K$ is a cycle,
all vertices of $K$ have degree at least two in $G$).  Furthermore, $K$ is an induced cycle, since $G\neq \Csixc$ is triangle-free and contains no dangerous cycles.
We can assume that $G$ is connected; otherwise, $G$ has a component $G_0$ disjoint from $K$, and this component either has a vertex of degree at most two forming \conf{1}
not interfering with the outer face, or a face bounded by a $(\le\!5)$-cycle $K_0$; in the latter case, we can consider $G_0$ drawn with $K_0$ as its outer face instead of $G$.

\begin{figure}
	\centering
        \subcaptionbox{$\Csixc$}[0.3\textwidth]{%
	\begin{tikzpicture}[scale=1.2]
				\draw (0,0) node[vertex] (v1) {}
		-- ++(0:1) node[vertex] (v2) {}
		-- ++(300:1) node[vertex] (v3) {}
		-- ++(240:1) node[vertex] (v4) {}
		-- ++(180:1) node[vertex] (v5) {}
		-- ++(-240:1) node[vertex] (v6) {}
		-- (v1);
		\draw (v3)-- ++ (v6);
	\end{tikzpicture}}
        \subcaptionbox{$\Csixv$}[0.3\textwidth]{%
	\begin{tikzpicture}[scale=1.2]
		\draw (0,0) node[vertex] (v1) {}
		-- ++(0:1) node[vertex] (v2) {}
		-- ++(300:1) node[vertex] (v3) {}
		-- ++(240:1) node[vertex] (v4) {}
		-- ++(180:1) node[vertex] (v5) {}
		-- ++(-240:1) node[vertex] (v6) {}
		-- (v1);

		\draw (v1)-- ++ (300:1) node[vertex,label=above:$v$] (v7) {} ;
		\draw (v3)-- ++ (v7);
		\draw (v5)-- ++ (v7);
	\end{tikzpicture}}
	\caption{Exceptional graphs in Lemma~\ref{lemma-unav}.}\label{fig-C6spec}
\end{figure}
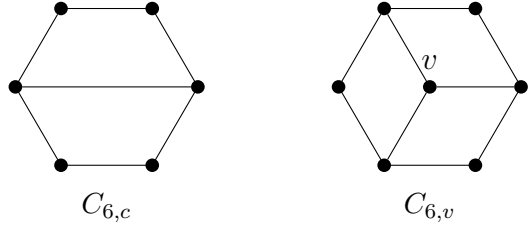

We now proceed by a discharging argument.  
Each vertex $v$ gets initial charge $c_0(v)=\deg(v)-4$, and each face $f$ gets initial charge $c_0(f)=|f|-4$.  By Euler's formula, the sum of the initial charges is
\begin{align*}
\sum_{v\in V} {(\deg(v)-4)}+\sum_{f\in F} {(|f|-4)}&=(2|E|-4|V|)+(2|E|-4|F|)\\
&=4(|E|-|V|-|F|)=-8.
\end{align*}
Next, we redistribute the charge according to the following rules.
\begin{enumerate}
	\item[\bf Rule 0:]
	        A non-outer face incident with a vertex $v\in V(K)$ of degree two sends $1/3$ to $v$.
	\item[\bf Rule 1:]
		Each (non-outer) face incident with an internal vertex $v$ of degree three sends $1/3$ to $v$.
	\item[\bf Rule 2:]
	        Let $f$ be a non-outer $4$-face incident with $k\ge 1$ vertices of $V(K)$, and let $v\in V(K)$ be a vertex incident with $f$.
		If $f$ is incident with an internal vertex of degree three, then $v$ sends $\frac{1}{3k}$ to $f$.
	\item[\bf Rule 3:]
	        Let $f$ be a (non-outer) $5$-face sharing an edge $uv$ with a $6$-face $g$, where $u$ and $v$ are internal vertices of degree $3$.
		Then $g$ sends $1/3$ to $f$.
	\item[\bf Rule 4:]
	        Let $f$ be a (non-outer) $5$-face, let $u$ be an internal vertex of degree three incident with $f$, and let $v$ be the neighbor of $u$
		not incident with $f$.  If $v\in V(K)$, then $v$ sends $1/3$ to $f$.
\end{enumerate}
Let $c$ denote the final charge obtained from $c_0$ by applying all the rules.  Note that no charge is created or lost,
and thus the sum of the final charges is still $-8$.  The charge of the outer face is unchanged, equal to $|K|-4$.

Let us first analyze the charge of a non-outer face $f=v_1\ldots v_{|f|}$.  If $|f|\ge 7$, then $f$ only sends charge by Rules 0 and 1 to incident vertices,
and thus $c(f)\ge c_0(f)-|f|/3=\frac{2}{3}|f|-4>0$.

Suppose that $|f|=6$.  If $f$ does not send charge by Rule 3, then $c(f)\ge \frac{2}{3}|f|-4=0$.  Let us consider the case that
$f$ sends charge by Rule 3 say to the face sharing the edge $v_2v_3$; hence, $v_2$ and $v_3$ are internal vertices of degree three.
If $v_1$ is an internal vertex, then $v_2$ and its three neighbors form \conf{2} ($G$ contains no path $P$ of length three
between $v_1$ and $v_3$, as $P$ together with the path $v_1v_2v_3$ would form a dangerous $5$-cycle) that does not interfere with the outer face.
Similarly, we can exclude the case that $v_4$ is internal.  Since $v_2$ and $v_3$ are internal vertices, it follows that $v_1$ and $v_4$ are not
vertices of $K$ of degree two.  Hence, $f$ does not send charge to $v_1$ and $v_4$, sends at most $1/3$ to each of $v_2$, $v_3$, $v_5$, and $v_6$,
and sends at most $2/3$ in total by Rule 3 to the faces incident with edges $v_2v_3$ and $v_5v_6$.  It follows that $c(f)\ge c_0(f)-2=0$.

Next, suppose that $|f|=5$.  Since $G\neq K$, observe that $f$ is incident with at most three vertices of $K$ of degree two,
and if $f$ is incident with at least one vertex of degree two, then it is incident with at least two vertices of $K$ of degree at least three.
In this case, $f$ sends $1/3$ to at most three vertices by Rules 0 and 1, and $c(f)\ge c_0(f)-1=0$.  Hence, we can assume that $f$ is
incident with no vertices of degree two.  Let $p$ be the number of internal vertices of degree three incident with $f$ whose neighbor
not incident with $f$ is internal, and let $q$ be the number of $6$-faces that share with $f$ an edge joining two internal vertices of degree three.
By Rules 1, 3, and 4, we have $c(f)\ge c_0(f)-(p-q)/3=(3+q-p)/3$, and thus if $c(f)<0$, then $p\ge q+4$; i.e., either $p=4$ and $q=0$, or $p=5$ and $q\le 1$.

Hence, we can assume that $v_1,\ldots, v_4$ are internal vertices of degree three such that their neighbors $u_1$, \ldots, $u_4$ not incident with $f$
are internal, and that the edge $v_2v_3$ is not incident with a $6$-face.
If $G-V(f)$ contains a path of length at most $2$ between $u_1$ and $u_4$,
then let $C$ be the $6$-cycle consisting of this path and the path $u_1v_1v_5v_4u_2$.  Since $C$ is not dangerous, the disk bounded by $C$ cannot contain
$f$, and thus $v_5$ is an internal vertex and $G_C$ contains all its neighbors.  Since $\deg(v_5)\ge 3$, $G_C$ is either $\Csixc$ or $\Csixv$, and in
either case, $v_5$ is an internal vertex of degree three and $v_1v_5$ is incident with a $4$-face.  However, this implies that $G$ contains \conf{5}
that does not interfere with the outer face.

Thus, we can assume that $G-V(f)$ contains no path of length at most $2$ between $u_1$ and $u_4$.  Similar argument shows that $u_1u_2, u_3u_4\not\in E(G)$
and (using the fact that $v_2v_3$ is not incident with a $6$-face) $G-V(f)$ contains no path of length $1$ or $3$ between $u_2$ and $u_3$.
Therefore, $f$ forms an appearance of \conf{4} in $G$, and since $u_1$, \ldots, $u_4$ are internal, this configuration does not interfere with the outer face.

Finally, suppose that $|f|=4$.  If say $v_1$ is a vertex of $K$ of degree two, then since $K$ is an induced cycle, it follows that $v_3$ is an internal
vertex.  Let $C$ be the cycle in $K+v_2v_3v_4$ distinct from $K$ and the boundary of $f$.  Since $C$ is not dangerous and $\deg(v_3)\ge 3$, we conclude
that $G_C$ is either $\Csixc$ or $\Csixv$.  The former is excluded, since $G$ is not isomorphic to $\Csixv$.  In the latter case, $G$ contains configuration \conf{5}
not interfering with the outer face.  This is a contradiction, and thus no vertex of $f$ has degree two, and $f$ sends no charge by Rule 0.
Let $p$ denote the number of internal vertices of degree three incident with $f$.  If $p\ge 2$, then $G$ contains \conf{3} or \conf{5} not interfering with the outer face.
If $p=0$, then $c(f)=c_0(f)=0$.  Hence, suppose that $p=1$, and say $v_1$ is an internal vertex of degree three.  Since $G$ does not contain dangerous cycles,
no path between $v_2$ and $v_4$ has length $3$.  Thus, $v_1$ and its neighbors form an appearance of \conf{2}, which must interfere with the outer face.
Consequently, at least one of $v_2$ or $v_4$ belongs to $V(K)$, and $c(f)=0$ by Rule 2.

\medskip

The preceding case analysis shows that the final charge of non-outer faces is non-negative. Let us now consider an internal vertex $v\in V(G)$.
If $\deg(v)\ge 4$, then $v$ neither sends nor receives charge and $c(v)=c_0(v)\ge 0$.  If $\deg(v)=3$, then $v$ receives
charge from all incident faces by Rule 1, and $c(v)=c_0(v)+1=0$.

Finally, let $v\in V(K)$ be a vertex incident with the outer face.  If $\deg(v)=2$, then $v$ receives $1/3$ by Rule 0 and does not send any charge (as we argued
before, non-outer 4-faces are not incident with degree 2 vertices, and thus Rule 2 does not apply), and $c(v)=c_0(v)+1/3=-5/3$.
If $\deg(v)\ge 3$, then $v$ sends at most $1/6$ to each of the $2$ incident non-outer faces sharing an edge with $K$ by Rule 2,
at most $1/3$ to each of the $\deg(v)-3$ other incident non-outer faces by Rule 2, and at most $1/3$ for each of $\deg(v)-2$ incident internal
vertices by Rule 4, giving the final charge $c(v)\ge c_0(v)-\frac{2}{3}(\deg(v)-2)=(\deg(v)-8)/3\ge -5/3$.

\medskip

In summary, all non-outer faces and internal vertices of $G$ have non-negative final charge and each vertex $v$ incident with the outer face
has final charge at least $-5/3$.  Furthermore, $c(v)=-5/3$ only if $\deg(v)=2$, or if $\deg(v)=3$, both incident non-outer faces have length $4$,
and $v$ is adjacent to an internal vertex of degree three.

It follows that the sum of the final charges is greater or equal
to the sum of the final charges of the outer face and its incident vertices, which is at least $|K|-4-\frac{5}{3}|K|=-4-\frac{2}{3}|K|$.
Since the sum of the final charges is $-8$, we conclude that $|K|=6$ and all vertices incident with $K$ have final charge $-5/3$.
This is only possible if every vertex $v\in V(K)$ has degree $2$ or $3$, all non-outer faces that share edge with $K$ have length $4$,
and all internal vertices with a neighbor in $K$ have degree three (not all vertices of $K$ have degree $2$ since $G$ is connected and $G\neq K$).
Since $K$ is an induced cycle and $G$ does not contain \conf{5} not interfering with the outer face, we conclude that each $4$-face whose boundary
intersects $K$ shares exactly $2$ edges with $K$, and thus $G$ is isomorphic to $\Csixv$.  This is a contradiction.
\end{proof}

\begin{corollary}\label{cor-unav}
Every plane triangle-free graph contains one of the configurations \conf{1}, \ldots, \conf{5}.
\end{corollary}
\begin{proof}
Let $G$ be a plane triangle-free graph, without loss of generality connected.  If $G$ contains a vertex of degree at most $2$, then \conf{1} appears in $G$.
Hence, assume that the minimum degree of $G$ is at least three.  Then, $G$ contains a face bounded by a cycle $K$ of length at most $5$.
Re-draw $G$ if necessary so that $K$ bounds the outer face of $G$.
Since the minimum degree of $G$ is at least three, $G$ is not a cycle, $\Csixc$ or $\Csixv$.

Let $K_1$ be a dangerous $(\le\!5)$-cycle in $G$ such that $G_{K_1}$ is minimal (we set $K_1=K$ if no $(\le\!5)$-cycle in $G$ is dangerous).
Let $K_2$ be a dangerous cycle of $G$ with $K_2\subseteq G_{K_1}$ such that the number of vertices of $G_{K_2}$ is minimum
(we set $K_2=K$ if $G$ contains no dangerous cycle).  By Lemma~\ref{lemma-unav}, $G_{K_2}$ contains one of the reducible configurations \conf{1}, \ldots, \conf{5}
that does not interfere with its outer face $K_2$.  Let $\gamma$ denote this configuration.  We claim that $\gamma$ is also a reducible configuration in $G$.
Let us discuss the configurations separately.

Suppose that $\gamma$ is \conf{2}; i.e., $G_{K_2}$ contains a vertex $v$ of degree three with neighbors $u$, $w$, $w'$ such that there exists no path of length $3$ between $w$ and $w'$ in $G_{K_2}$,
and since $\gamma$ does not interfere with the outer face of $G_{K_2}$, we have $v,w,w'\not\in V(K_2)$.  Hence, $v$ has degree $3$ in $G$ as well.  Furthermore, if there exists a path of length $3$
between $w$ and $w'$ in $G$, then there exist adjacent vertices $z,z'\in V(K_2)$ such that $zz'\not\in E(G_{K_2})$ and $wz,w'z'\in E(G)$.  Since $zz'$ is a chord of $K_2$ and $G$ is triangle-free,
it follows that $|K_2|=6$ and $z$ and $z'$ are opposite vertices of $K_2$, i.e., $K_2=zx_1x_2z'y_1y_2$.  Since $|K_1|\le 5$, the cycle $K_1$ has no chord, and since $z,z'\in V(G_{K_2})\subseteq V(G_{K_1})$,
we conclude that $zz'\in E(G_{K_1})$.  Since the $5$-cycle $C=vwzz'w'$ separates $x_1,x_2$ from $y_1,y_2$, it is dangerous and $G_C\subsetneq G_{K_1}$; this contradicts the choice of $K_1$.
Therefore, $G$ contains no path of length $3$ between $w$ and $w'$, and thus $\gamma$ also forms \conf{2} in $G$.

Suppose that $\gamma$ is \conf{4}; i.e., $G_{K_2}$ contains a $5$-face $f=v_1\ldots v_5$ with $v_1, \ldots, v_4$ having degree three such that,
denoting for $i=1,\ldots, 4$ the neighbor of $v_i$ outside of $f$ by $u_i$, $G_{K_2}-V(f)$ contains no path of length at most $2$ between $u_1$ and $u_4$,
and no path of length $1$ or $3$ between $u_2$ and $u_3$, and $u_1u_2,u_3u_4\not\in E(G_{K_2})$, and $v_1,\ldots, v_4,u_1,\ldots, u_4\not\in V(K_2)$.
Clearly, $\gamma$ forms \conf{4} in $G$ as well, unless $G-V(f)$ contains a path $u_2zz'u_3$ of length three.  As in the previous paragraph,
this is only possible if $K_2=zx_1x_2z'y_1y_2$, and letting $C$ be the $6$-cycle $u_3v_3v_2u_2zz'$, we have $G_C\subseteq G_{K_1}$, and since $C$ separates $x_1,x_2$ from $y_1,y_2$, it is dangerous.
Considering the $4$-cycles $C_1=zx_1x_2z'$ and $C_2=zy_1y_2z'$, we can by symmetry assume that the closed disk bounded by $C_2$ contains both $C$ and $C_1$ (and by the minimality in the choice of $K_1$,
we have $C_2=K_1$).  By the minimality in the choice of $K_1$, the cycle $C_1$ is not dangerous, and thus $C_1$ bounds a face.  Therefore, $V(G_C)\subseteq V(G_{K_2})\setminus\{y_1,y_2\}$, which contradicts
the minimality in the choice of $K_2$.

Finally, if $\gamma$ is \conf{1}, \conf{3}, or \conf{5}, then the vertices of $\gamma$ whose degree is required to be equal to $2$ or $3$ are not incident with $K_2$, and thus their
degree in $G_{K_2}$ is the same as their degree in $K$.  Consequently, $\gamma$ is a reducible configuration in $G$ as well.
\end{proof}

\section{Independent sets}

Our main result is now an easy consequence.

\begin{proof}[Proof of Theorem~\ref{thm-main}]
Suppose for a contradiction that the claim is false, and there exists a tight graph $G\not\in \G$.  Choose such a graph with the minimum number of vertices, so that $G$ is a minimum counterexample.
By Corollary~\ref{cor-unav}, $G$ contains one of the configurations \conf{1}, \ldots, \conf{5}, which contradicts Lemma~\ref{lemma-noredu}.
\end{proof}

\bibliographystyle{siam}
\bibliography{indep}

\end{document}